# Pareto front generation with knee-point based pruning for mixed discrete multi-objective optimization

Juseong Lee[1], Sang-Il Lee[2], Jaemyung Ahn[1*], and Han-Lim Choi[1]

1: *Korea Advanced Institute of Science and Technology* (*KAIST*), *Daejeon, Republic of Korea*, 2: *Korea Aerospace Research Institute* (*KARI*), *Daejeon, Republic of Korea*, *: Corresponding Author*

**Abstract** This note proposes an algorithm to generate the Pareto front of a mixed discrete multi-objective optimization problem based on the pruning of irrelevant subproblems. An existing pruning-based method for a mixed discrete bi-objective problem is extended for general multi-objective cases by introducing a new reference point for pruning decision – the knee point. The validity of the proposed procedure is demonstrated through case studies.

## 1. Introduction

This note concerns a mixed-discrete multi-objective optimization (MOO) problem defined as follows:

$$\min_{\mathbf{x}} \mathbf{J}(\mathbf{x}) = \min_{\mathbf{y},\mathbf{z}} \mathbf{J}(\mathbf{y},\mathbf{z}) = [J_1(\mathbf{y},\mathbf{z}),\cdots,J_m(\mathbf{y},\mathbf{z})]^T \quad \textbf{(P)}$$

subject to

$$\mathbf{g}(\mathbf{y},\mathbf{z}) \leq \mathbf{0}, \quad \mathbf{h}(\mathbf{y},\mathbf{z}) = \mathbf{0}$$
$$l_i \leq y_i \leq u_i \quad (i=1,\cdots,n_c)$$
$$z_j \in Z_j = \{z_{j,1},\cdots,z_{j,d_j}\} \quad (j=1,\cdots,n_d)$$

where $\mathbf{x} = [\mathbf{y}\ \mathbf{z}]$ is the design vector, $\mathbf{y}/\mathbf{z}$ are the continuous/discrete components of $\mathbf{x}$ whose dimensions are respectively $n_c/n_d$, $m$ is the number of objectives, $\mathbf{g}/\mathbf{h}$ are the inequality/equality constraint vectors, and $Z_j$ is the set of values that $j^{th}$ discrete design variable ($z_j$) can take ($|Z_j| = d_j$). The Pareto optimal solution of the problem ($\mathcal{X}^*$) is defined as follows:

$$\mathcal{X}^* = \{\mathbf{x}^* \in \mathcal{X} \mid \nexists \mathbf{x} \in \mathcal{X}\setminus\{\mathbf{x}^*\}\ s.t.\ J_i(\mathbf{x}) \leq J_i(\mathbf{x}^*)\ \forall i,\ J_j(\mathbf{x}) < J_j(\mathbf{x}^*)\ \exists j\} \quad (1)$$

where $\mathcal{X}$ denotes the set of feasible design vectors. In addition, the Pareto front ($\mathcal{J}^*$) is defined as the points in objective space corresponding to $\mathcal{X}^*$ as follows:

$$\mathcal{J}^* = \{\mathbf{J}(\mathbf{x}) \mid \mathbf{x} \in \mathcal{X}^*\} \quad (2)$$

Hong et al. (2015) proposed a procedure to solve the bi-objective version of this problem ($m = 2$) based on two-phase pruning to eliminate irrelevant subproblems. A subproblem instantiated by specifying the discrete design variable is pruned out if its reference point (Phase A: utopia point, Phase B: center point) is dominated by the "master front," which is defined in Section 2. Their method performed effectively for mixed-discrete bi-objective optimization (BOO) problems. However, their formulation could handle only two objectives, and its generalization for three or more objectives was not addressed in the paper. In addition, the center point used in their second phase may cause unsuccessful pruning decision depending on the shape of the resultant Pareto front.

This paper proposes formulation/procedure that can handle three or more objective functions and utilizes a new reference point (the knee point) for Phase B pruning to overcome the aforementioned disadvantages of the previous study.

## 2. Review: Pareto-front generation for mixed discrete bi-objective optimization with center-point based pruning

The Pareto front generation procedure with center point based pruning by Hong et al. (2015) – referred to as *C-Pruning* in the rest of this note – is composed of two phases that adopt the utopia points and the center points as the references for pruning decision, respectively.

[Phase A of *C-Pruning*]

(Step A-1) *Computing anchor/utopia points of subproblems*: The *C-Pruning* first generates $\mathcal{Z}$, the set of all feasible discrete variable combinations, and its index set $\mathcal{K}$. Then the subproblem $\mathbf{P}_k$ is instantiated by setting $\mathbf{z} = \mathbf{z}_k \in \mathcal{Z}$ in **P** as follows:

$$\min_{\mathbf{y}} \mathbf{J}(\mathbf{y},\mathbf{z}_k) = [J_1(\mathbf{y},\mathbf{z}_k),\cdots,J_m(\mathbf{y},\mathbf{z}_k)]^T \quad \textbf{(P}_k\textbf{)}$$

subject to

$$\mathbf{g}(\mathbf{y},\mathbf{z}_k) \leq \mathbf{0}, \quad \mathbf{h}(\mathbf{y},\mathbf{z}_k) = \mathbf{0}$$
$$l_i \leq y_i \leq u_i \quad (i=1,\cdots,n_c)$$

The anchor points of $\mathbf{P}_k$ are defined as follows:

$$\mathbf{J}_k^{A,i} = \mathbf{J}(\mathbf{y}_k^{A,i},\mathbf{z}_k) \equiv [J_{k,1}^{A,i},\cdots,J_{k,j}^{A,i},\cdots,J_{k,m}^{A,i}]^T \quad (i=1,\cdots,m) \quad (3)$$

where $\mathbf{y}_k^{A,i} = \arg\min_{\mathbf{y}\in\mathcal{Y}_k} J_i(\mathbf{y},\mathbf{z}_k)$, or the design vector that minimizes the $i^{th}$ objective. The utopia point of $\mathbf{P}_k$ is expressed as

$$\mathbf{J}_k^U = [J_{k,1}^{A,1}, J_{k,2}^{A,2}, \cdots, J_{k,m}^{A,m}]^T \qquad (4)$$

(Step A-2) *Generating a master front*: The index set of subproblems with efficient (non-dominated) utopia points ($\mathcal{K}_1^M$) is identified using pairwise comparisons of utopia points as follows:

$$\mathcal{K}_1^M = \{k \mid \nexists l \in \mathcal{K}\setminus\{k\} \text{ s.t. } \mathbf{J}_l^U \leq \mathbf{J}_k^U\} \qquad (5)$$

The solution of $\mathbf{P}_k$ ($\mathcal{Y}_k^*$) and associated sub Pareto front ($\mathcal{J}_k^*$) are defined as follows:

$$\mathcal{Y}_k^* = \{\mathbf{y}^* \in \mathcal{Y}_k \mid \nexists y \in \mathcal{Y}_k \setminus \{\mathbf{y}^*\} \text{ s.t.}$$
$$J_i(\mathbf{y}, \mathbf{z}_k) \leq J_i(\mathbf{y}^*, \mathbf{z}_k) \ \forall i, \qquad (6)$$
$$J_j(\mathbf{y}, \mathbf{z}_k) < J_j(\mathbf{y}^*, \mathbf{z}_k) \ \exists j\}$$

$$\mathcal{J}_k^* = \{\mathbf{J}(\mathbf{y}, \mathbf{z}_k) \mid y \in \mathcal{Y}_k^*\} \qquad (7)$$

where $\mathcal{Y}_k$ denotes the set of feasible design vectors of $\mathbf{P}_k$. In addition, we define a master front ($\mathcal{J}^M$) as the set of non-dominated solutions out of the collection of sub Pareto fronts ($\mathcal{J}_k^*$) for $k \in \mathcal{K}_1^M$.

The proposed algorithm (*K-Pruning*) assumes that the Pareto front of subproblem ($\mathbf{P}_k$) can be found with NBI (Das and Dennis 1998; PK Shukla 2007) or AWS (Kim and de Weck 2006; Hwang and Masud 2012). This implies that a gradient-based method is used to solve the NLPs instantiated during NBI/AWS, which works only if the subproblem ($\mathbf{P}_k$) is differentiable. Note that the overall Pareto is can be non-differentiable or even discontinuous though because it is created by "patching" sub Pareto fronts ($\mathcal{J}_k^*$).

(Step A-3) *Utopia point based pruning*: Subproblem $\mathbf{P}_k$ is pruned out if its utopia point is dominated by the master front. The index set of pruned subproblems during the utopia point based pruning ($\mathcal{K}_\varnothing^U$) is defined as follows:

$$\mathcal{K}_\varnothing^U = \{k \mid \exists \mathbf{J}(\mathbf{x}) \in \mathcal{J}^M \text{ s.t. } \mathbf{J}(\mathbf{x}) \leq \mathbf{J}_k^U\} \qquad (8)$$

[Phase B of *C-Pruning*]

(Step B-1) *Computing center points of subproblems*: The center point of a subproblem ($\mathbf{J}_k^C = \mathbf{J}(\mathbf{y}_k^C, \mathbf{z}_k)$) is obtained by solving the following NLP.

$$\mathbf{y}_k^C = \arg\min_\mathbf{y} J^1(\mathbf{y}, \mathbf{z}_k) \qquad (9)$$

subject to constraints for $\mathbf{P}_k$ and

$$\frac{J_2(\mathbf{y}, \mathbf{z}_k) - J_{k,2}^{A,2}}{J_{k,2}^{A,1} - J_{k,2}^{A,2}} - \frac{J_1(\mathbf{y}, \mathbf{z}_k) - J_{k,1}^{A,1}}{J_{k,1}^{A,2} - J_{k,1}^{A,1}} \leq 0 \qquad (10)$$

Fig. 1-(a) shows geometric interpretation of the center point. The constraint representing the feasible region (Eq. (10)) is expressed as gray area. The center point can be found by minimizing $J_1$ within the feasible region. This algorithm is applicable to bi-objective case and the value of $m$ introduced in $\mathbf{P}_k$ is selected as 2.

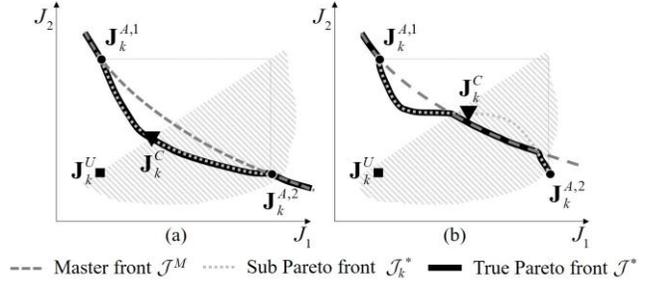

Figure 1: (a) Center point geometry, (b) Unsuccessful pruning with center point

(Step B-2) *Center point based pruning*: A subproblem whose center point ($\mathbf{J}_k^C$) is dominated by the master front ($\mathcal{J}^M$) is pruned out. The index set of eliminated subproblems during the center point based pruning ($\mathcal{K}_\varnothing^C$) is defined as

$$\mathcal{K}_\varnothing^C = \{k \mid \exists \mathbf{J}(\mathbf{x}) \in \mathcal{J}^M \text{ s.t. } \mathbf{J}(\mathbf{x}) \leq \mathbf{J}_k^C\} \qquad (11)$$

(Step B-3) *Generating Pareto front*: The (approximate) solution of the original problem is obtained by combining sub Pareto fronts for $k \in \mathcal{K}\setminus(\mathcal{K}_\varnothing^U \cup \mathcal{K}_\varnothing^C)$. Note this solution is exact if $\mathcal{K}_\varnothing^U \cup \mathcal{K}_\varnothing^C$ is identical to the set of all irrelevant subproblems ($\mathcal{K}_\varnothing$).

While the *C-Pruning* effectively prunes out irrelevant subproblems for practical mixed-discrete BOOs, it has some disadvantages as well. First, their formulation was developed in the context of two objectives. Its extension to three or higher dimension, which is not very straightforward, was not provided in their work. Secondly, the center point can omit parts of the true Pareto front depending on the shape of the true Pareto front, which case is illustrated in Fig. 1. *C-Pruning* considers $\mathcal{J}_k^*$ presented in Fig. 1-(b) as irrelevant because $\mathbf{J}_k^C$ is dominated by $\mathcal{J}^M$. However the bulged part of $\mathcal{J}_k^*$ is actually Pareto-optimal.

## 3. Knee point based pruning for Pareto front generation of a mixed-discrete multi-objective optimization

The knee point based pruning algorithm, which is referred to as *K-Pruning* in the rest of this note, modifies Phase B of the *C-Pruning* algorithm by changing its reference point for pruning decision (from center point to knee point). The knee point is conceptually the most bulged region, which has a couple of different definitions (Das 1999; Sudeng and Wattanapongsakorn 2015; Rachmawati and Srinivasan 2009).

[Phase B of *K-Pruning*]

(Step B-1) *Computing Knee points of subproblems*: The second phase of *K-Pruning* starts with solving the following optimization ($\mathbf{P}_k^K$) for $k \in \mathcal{K} \setminus (\mathcal{K}_1^M \cup \mathcal{K}_\emptyset^U)$, which is the distance based method suggested by Das (1999).

$$[\mathbf{y}_k^K, \boldsymbol{\beta}_k^K, t] = \arg\max_{\mathbf{y}, \boldsymbol{\beta}, t} t \quad (\mathbf{P}_k^K)$$

subject to constraints of $\mathbf{P}_k$ and additional constraints:

$$\boldsymbol{\Phi}_k \boldsymbol{\beta} + t\mathbf{n}_k = \mathbf{J}(\mathbf{y}, \mathbf{z}_k)$$

$$\sum_{j=1}^{m} \beta_j = 1, \quad \beta_j \geq 0 \quad (j = 1, \cdots, m)$$

where $\boldsymbol{\Phi}_k = [\mathbf{J}_k^{A,1}, \cdots, \mathbf{J}_k^{A,m}]$ is the payoff matrix of $\mathbf{P}_k$, $\boldsymbol{\beta} = [\beta_1, \ldots, \beta_m]$, and vector $\mathbf{n}_k$ is normal to plane $\boldsymbol{\Phi}_k \boldsymbol{\beta}$. Das et al. (1998) mentioned that if $\mathbf{n}_k$ is not available, a quasi-normal vector ($\tilde{\mathbf{n}}$) can be used for NBI method. For example, a vector normal to the plane $\boldsymbol{\Phi}\boldsymbol{\beta}$ can be suggested as $\tilde{\mathbf{n}}$, where $\boldsymbol{\Phi} = [\mathbf{J}^{A,1}, \ldots, \mathbf{J}^{A,m}]$ is the payoff matrix of $\mathbf{P}$ ($\mathbf{J}^{A,i}$: anchor point of $\mathbf{P}$ associated with $i^{th}$ objective).

The knee point of $\mathbf{P}_k$ ($\mathbf{J}_k^K$) is defined using $\mathbf{y}_k^K$ as

$$\mathbf{J}_k^K = \mathbf{J}(\mathbf{y}_k^K, \mathbf{z}_k) \quad (12)$$

Geometric interpretation of knee point is presented in Fig. 2.

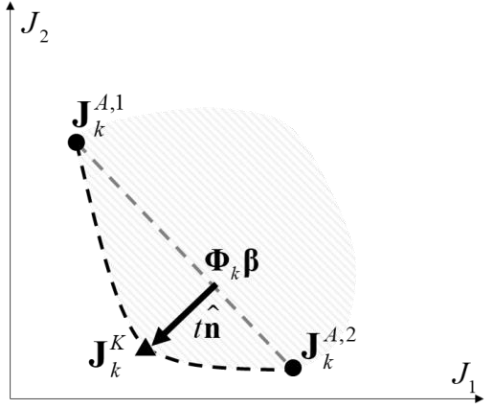

Figure 2: Geometric interpretation of a knee

(Step B-2) *Knee-point based pruning*: If $\mathbf{J}_k^K$ is dominated by $\mathcal{J}^M$, we can conclude that $\mathbf{P}_k$ is irrelevant.

$$\mathcal{K}_\emptyset^K = \{k \mid \exists \mathbf{J}(\mathbf{x}) \in \mathcal{J}^M \text{ s.t. } \mathbf{J}(\mathbf{x}) \leq \mathbf{J}_k^K\} \quad (13)$$

(Step B-3) *Generating Pareto front*: $\mathcal{J}^*$ can be generated by combining $\mathcal{J}_k^*$ for $k \in \mathcal{K} \setminus (\mathcal{K}_\emptyset^U \cup \mathcal{K}_\emptyset^K)$ and extracting non-dominated solutions.

The *K-Pruning* has three meaningful advantages over the *C-Pruning*. First, *K-Pruning* is applicable to problems with more than two objectives since the knee point is properly defined in general *m*-dimensional objective space. Secondly, the knee point provides an intuitive and geometrically easy-to-understand reference point for pruning decision even for a complicated sub Pareto front. Fig. 3 compares the reference points for pruning decision based on the center point and the knee point for three different sub Pareto front geometries – (a) convex and symmetric, (b) convex and skewed, and (c) nonconvex. Both methods provide proper references in Fig. 3-(a). However, in Figs 3-(b) and 3-(c), the center point does not represent the sub Pareto front for pruning decision as the knee point, which can lead to over-pruning of relevant subproblems. This issue will be discussed in the first case study (TP1). Lastly, *K-Pruning* can find all knee points of the original problem $\mathbf{P}$ if each subproblem is convex. Note that the knee regions are where the maximum trade-off of objective functions takes place and thus can provide attractive design alternatives (Branke et al. 2004; Deb and Gupta 2010).

The effectiveness of the proposed *K-Pruning* algorithm is demonstrated through case studies in the next section.

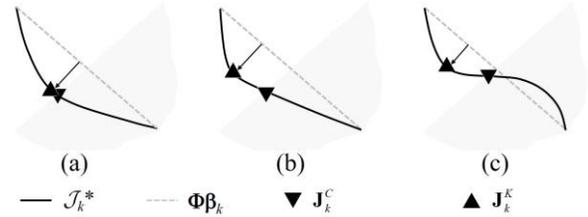

Figure 3: Reference point comparisons for two pruning methods

### 4. Case study
#### 4.1. Test Problem 1: Nonconvex sub Pareto fronts

Test problem 1 (TP1) is a bi-objective problem that can compare the performance of the two pruning methods. TP1 involves skewed/nonconvex sub Pareto fronts as presented in Figs. 3-(b) and 3-(c). It is mathematically formulated as:

$$\min_{y, \mathbf{z}} \mathbf{J} = \min_{y, \mathbf{z}} \begin{bmatrix} J_1(y, \mathbf{z}) \\ J_2(y, \mathbf{z}) \end{bmatrix} \quad (\mathbf{TP1})$$

$$-2 \leq y \leq 2, \quad z_1 \in \{1, 2, 3\}, \quad z_2, z_3 \in \{-1, 0, 1\}$$

where

$$J_1(y, \mathbf{z}) = C_{z_1} + z_2 - \tfrac{1}{2} z_3 + A_{z_1} y^3 + B_{z_1} y^2 + (1 - 4A_{z_1}) y$$

$$J_2(y, \mathbf{z}) = C_{z_1} - \tfrac{1}{2} z_2 - z_3 + A_{z_1} y^3 + B_{z_1} y^2 + (1 - 4A_{z_1}) y$$

$$[A_{z_1}, B_{z_1}, C_{z_1}] = \left[\tfrac{4}{45}, \tfrac{2}{45}, \tfrac{-6}{45}\right], \left[0, \tfrac{1}{20}, \tfrac{-1}{5}\right], \left[-\tfrac{4}{45}, \tfrac{2}{45}, \tfrac{-6}{45}\right] (z_1 = 1, 2, 3)$$

Note that $z_1$ is a categorical variable that collectively changes the coefficients of the objective functions. Its total number of subproblems ($|\mathcal{K}|$) is 27. Table 3 indicates 6 out of 22 subproblems eliminated in Phase B of *C-Pruning* were actually relevant (over-pruning), and 2 irrelevant problems were not pruned (under-pruning). On the other hand, *K-Pruning* did not over-prune any relevant subproblems while 2 irrelevant subproblems were under-pruned. The obtained Pareto front is presented in Fig. 5-(a).

Table 1: Pruning characteristics of two methods (TP1)

| Method | $|\mathcal{K}|$ | $|\mathcal{K}_1^m|$ | $|\mathcal{K}_\varnothing^U|$ | $|\mathcal{K}_\varnothing^K|$ | $|\mathcal{K}_\varnothing|$ | Over-pruning | Under-pruning |
|---|---|---|---|---|---|---|---|
| C-Pruning | 27 | 5 | 0 | 22 | 18 | 6 | 2 |
| K-Pruning | | | | 16 | | 0 | 2 |

### 4.2. Test Problem 2: DEB3DK with discrete variables

A tri-objective optimization problem referred to as DEB3DK (Branke et al. 2004) is modified and used as the second test problem (TP2). This problem was originally developed as a test case for studies to find knee regions (Rachmawati et al. 2004; and Bechikh et al. 2011). Note that the original formulation of DEB3DK involved only continuous variables, and some of its parameters are switched to discrete design variables (z) in this test problem. The problem is mathematically formulated as follows:

$$\min_{\mathbf{y},\mathbf{z}} \mathbf{J} = \min_{\mathbf{y},\mathbf{z}} \begin{bmatrix} J_1(\mathbf{y},\mathbf{z}) & J_2(\mathbf{y},\mathbf{z}) & J_3(\mathbf{y},\mathbf{z}) \end{bmatrix}^T \quad (\textbf{TP2})$$

subject to

$$J_1(\mathbf{y},\mathbf{z}) = g(\mathbf{y},\mathbf{z})r_1(\mathbf{y},\mathbf{z})\left(\sin\left(\frac{y_1\pi}{2}\right) + \cos\left(\frac{y_2\pi}{2}\right) + z_5 z_6\right)$$

$$J_2(\mathbf{y},\mathbf{z}) = g(\mathbf{y},\mathbf{z})r_2(\mathbf{y},\mathbf{z})\left(\sin\left(\frac{y_1\pi}{2}\right) + \cos\left(\frac{y_2\pi}{2}\right) + z_6 z_7\right)$$

$$J_3(\mathbf{y},\mathbf{z}) = g(\mathbf{y},\mathbf{z})r_3(\mathbf{y},\mathbf{z})\left(\cos\left(\frac{y_1\pi}{2}\right) + z_5 z_7\right)$$

$$g(\mathbf{y},\mathbf{z}) = 1 + \frac{9}{n_c - 1}\sum_{i=2}^{n_c} y_i$$

$$r_1(\mathbf{y},\mathbf{z}) = 5 + z_4(y_1 - z_1)^2 + \frac{\cos(2z_0 \pi y_1)}{z_0}$$

$$r_2(\mathbf{y},\mathbf{z}) = 5 + (20 - z_4)(y_2 - z_2)^2 + \frac{\cos(2z_0 \pi y_2)}{z_0}$$

$$r(\mathbf{y},\mathbf{z}) = r_1 + r_2(1 - z_2)$$

$$0 \leq y_i \leq 1, \quad \forall i \in \{1, \cdots, 30\}$$

$$z_i \in Z_i \quad \text{for} \quad \forall i \in \{0, \cdots, 7\}$$

$$Z_0 = \{1, 2\}, Z_1, Z_2, Z_3 = \{0.2, 0.5, 0.8\}, Z_4 = \{5, 10, 15\}$$

$$Z_5, Z_6, Z_7 = \{0, 1, 2, 3\}$$

Pruning characteristics of TP2 is summarized in Table 2. TP2 has total 34,992 ($2 \times 3^4 \times 6^3$) and K-Pruning concluded that only 2,838 subproblems are relevant, which accounts for only 8% of the total number. Fig.5-(b) shows the Pareto front generated – 85 subproblems are relevant, and K-Pruning found all relevant subproblems.

Table 2: Pruning characteristics of TP2

| $|\mathcal{K}|$ | $|\mathcal{K}_1^M|$ | $|\mathcal{K}_\varnothing^U|$ | $|\mathcal{K}_\varnothing^K|$ | $|\mathcal{K}_\varnothing|$ | Over-pruning | Under-Pruning |
|---|---|---|---|---|---|---|
| 34,992 | 2,592 | 24,948 | 7,206 | 32,154 | 0 | 2,753 |

### 4.3. Test Problem 3: Tri-Objective Nine Bar Truss Problem

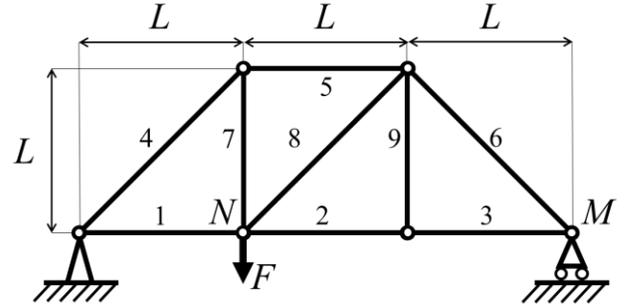

Figure 4: Structural configuration of nine bar truss (TP2)

The third test problem (TP3) is created by modifying a "classic" structural optimization problem called *Nine Bar Truss Problem* (e.g. Mela et al. 2007; Hong et al. 2014). The structural configuration of the nine bar truss is presented in Fig.4. Cross sectional areas of trusses 1-3 ($x_1 - x_3$), and discrete choices of cross sectional area / material type combination for trusses 4-9 ($x_4 - x_9$) are defined as decision variables of the problem. Note that variables $x_1$, $x_2$, and $x_3$ are continuous and the other variables are discrete ($\mathbf{y} = [x_1\ x_2\ x_3]$, $\mathbf{z} = [x_4\ x_5\ x_6\ x_7\ x_8\ x_9]$). The problem is formulated using parameters summarized in Tables 3 and 4 as follows:

$$\min_{\mathbf{x}} \mathbf{J} = \min \begin{bmatrix} J_1(\mathbf{x}) \\ J_2(\mathbf{x}) \\ J_3(\mathbf{x}) \end{bmatrix} = \begin{bmatrix} L\left(\sum_{i=1}^9 l_i C(x_i)\right) \\ \frac{FL}{9}\left(\sum_{i=1}^9 \frac{b_i}{A(x_i)E(x_i)}\right) \\ \frac{FL}{3}\left(\frac{2}{x_1 E_0} + \frac{1}{x_2 E_0} + \frac{1}{x_3 E_0}\right) \end{bmatrix} \quad (\textbf{TP3})$$

subject to

$$\tfrac{2}{3} \leq x_1 \leq 10, \quad \tfrac{1}{3} \leq x_2 \leq 10, \quad \tfrac{1}{3} \leq x_3 \leq 10$$

$$x_i \in \{1,2,3,4,5\} \quad (i = 4,\cdots,9)$$

Three objectives to be minimized are the following: *material cost* ($J_1$) determined by length ($l_i$) and price per length ($C(x_i)$), *vertical displacement of node N* ($J_2$), and the *horizontal displacement of node M* ($J_3$). The cross-section area, Young's modulus ($E(x_i)$), and unit cost ($C(x_i)$) of truss $i$ are given as follows:

$$A(x_i) = \begin{cases} x_i, & 1 \leq i \leq 3 \\ A_{x_i}, & 4 \leq i \leq 9 \end{cases} \qquad E(x_i) = \begin{cases} E_0, & 1 \leq i \leq 3 \\ E_{x_i}, & 4 \leq i \leq 9 \end{cases}$$

$$C(x_i) = \begin{cases} x_i C_0, & 1 \leq i \leq 3 \\ C_{x_i}, & 4 \leq i \leq 9 \end{cases}$$

Table 3: Coefficients in TP3

| $i$ | 1 | 2 | 3 | 4 | 5 | 6 | 7 | 8 | 9 |
|---|---|---|---|---|---|---|---|---|---|
| $l_i$ | 1 | 1 | 1 | $\sqrt{2}$ | 1 | $\sqrt{2}$ | 1 | $\sqrt{2}$ | 1 |
| $b_i$ | 4 | 1 | 1 | $8\sqrt{2}$ | 4 | $2\sqrt{2}$ | 4 | $2\sqrt{2}$ | 0 |

Table 4: Parameter values used in TP3

| Identifier ($z$) | Cross-section area ($A_z$, in$^2$) | Young's Modulus ($E_z$, ksi) | Unit cost ($C_z$, $/ft) |
|---|---|---|---|
| 0 | N/A | 45 | 6.00 |
| 1 | 1.00 | 29 | 3.12 |
| 2 | 0.25 | 29 | 0.85 |
| 3 | 1.00 | 45 | 5.70 |
| 4 | 0.39 | 45 | 3.19 |
| 5 | 0.25 | 64 | 2.59 |

The number of unique discrete variable combinations ($|\mathcal{K}|$) is 15,625 ($6^5$). The proposed K-Pruning procedure deleted 5,534 subproblems – including 5 over-pruning and 52 under-pruning. The obtained Pareto front and over-pruned solution points are presented in Fig. 5-(c). The pruning characteristics of K-pruning for TP3 is presented in Table 5.

Note that $J_2$ and $J_3$ of TP3 are not in significant conflict. This is reflected in the projection of the Pareto front in $J_2$-$J_3$ plane, which does not present clear trade-off relationship. In practice, it is desirable that the trade-off between objectives of MOO are clearly explainable to obtain meaningful study results.

Table 3: Pruning characteristics TP3

| $|\mathcal{K}|$ | $|\mathcal{K}_1^m|$ | $|\mathcal{K}_\varnothing^U|$ | $|\mathcal{K}_\varnothing^K|$ | $|\mathcal{K}_\varnothing|$ | Over-pruning | Under-Pruning |
|---|---|---|---|---|---|---|
| 15,625 | 38 | 0 | 11,835 | 11,835 | 0 | 3,752 |

## 5. Discussion: Computational Cost and Limitation

The primary target of the proposed approach is the problems including "categorical" discrete variables, which cannot be relaxed as real numbers (e.g. material types: steel/aluminum/tungsten, moving mechanism: wheel/legs). While there are efficient algorithms based on relaxation of discrete variables, the algorithms do not work on problems involving the categorical variables because they cannot be relaxed.

It should be noted that the increase in the numbers of discrete variables/options that each of them can take can result in the explosion in prohibitively large number of subproblems that we have to solve. For example, if there are 10 discrete variables each of which can take 10 values, we have to find sub Pareto fronts associated with $10^{10}$ subproblems. Such a case can be addressed using multi-objective heuristics (e.g. multi-objective genetic algorithm); its performance would be highly dependent on the nature of the discrete variables as well.

## 6. Conclusion

A procedure to generate the Pareto front of mixed-discrete multi-objective optimization (MOO) problems using the knee-point based pruning (*K-Pruning*) is proposed. The procedure is developed by extension of existing method to solve a mixed-discrete bi-objective optimization (BOO) using the center-point based pruning (*C-Pruning*). In addition to the capability to handle three or more objectives, the proposed method provides more intuitive reference points for pruning decision than the existing method, and makes sure that the region where high degree of trade-off between different objectives is included in the solution. Case study with three test problems demonstrated the effectiveness of the proposed procedure for solving the multi-objective MOO.

## References

Branke J, Deb K, Dierolf H, Osswald M (2004) Finding knees in multi-objective optimization. In: Parallel Problem Solving from Nature (PPSN) 2004 - Lecture Notes in Computer Science Vol 3242: 722–731

Bechikh, S., Said, L. Ben, & Ghédira, K. (2011). Searching for knee regions of the Pareto front using mobile reference points. Soft Computing, 15(9), 1807–1823. https://doi.org/10.1007/s00500-011-0694-3

Das I (1999) On characterizing the 'knee' of the Pareto curve based on normal-boundary intersection. Structural Optimization 18 (2):107-115

Das I, Dennis JE (1998) Normal-boundary intersection: a new method for generating the Pareto surface in nonlinear multi-criteria optimization problems. Society for Industrial and Applied Mathematics Journal on Optimization 8(3): 631–657

Deb K, Gupta S (2010) Towards a link between knee Solutions and preferred solution methodologies. In: Panigrahi B, Das S, Suganthan P, Dash S (eds) Swarm, Evolutionary, and Memetic Computing (SEMCCO) 2010 – Lecture Notes in Computer Science Vol 6466:182-189

Hong S, Ahn J, Choi HL (2014) Pruning-based Pareto front generation for mixed-discrete bi-objective optimization. Structural and Multidisciplinary Optimization 51(1):193–98

Hwang Cl, Masud A (2012) Multiple objective decision making - methods and applications: a state-of-the-art survey. Springer

Kim IY and de Weck O (2006) Adaptive weighted sum method for multiobjective optimization: a new method for Pareto front generation. Structural and Multidisciplinary Optimization 31(2):105-116

Mela K, Koski J, Silvennoinen R (2007) Algorithm for generating the Pareto optimal set of multiobjective nonlinear mixed-integer optimization problems. In: 48th AIAA/ASME/ASCE/AHS/ASC Structures, Structural Dynamics, and Materials Conference. Structures, Structural Dynamics, and Materials and Co-Located Conferences. AIAA, Honolulu, Hawaii

Rachmawati L, Srinivasan D (2009) Multiobjective evolutionary algorithm with controllable focus on the knees of the Pareto front. IEEE Transactions on Evolutionary Computation 13(4):810–24

Shukla, P. (2007). On the normal boundary intersection method for generation of efficient front. Computational Science – ICCS 2007. Lecture Notes in Computer Science, 4487. Springer, Berlin, Heidelberg.

Sudeng S, Wattanapongsakorn N (2015) Finding knee solutions in multi-objective optimization using extended angle dominance approach. In: Information Science and Applications. Lecture Notes in Electrical Engineering Vol 339:673-679

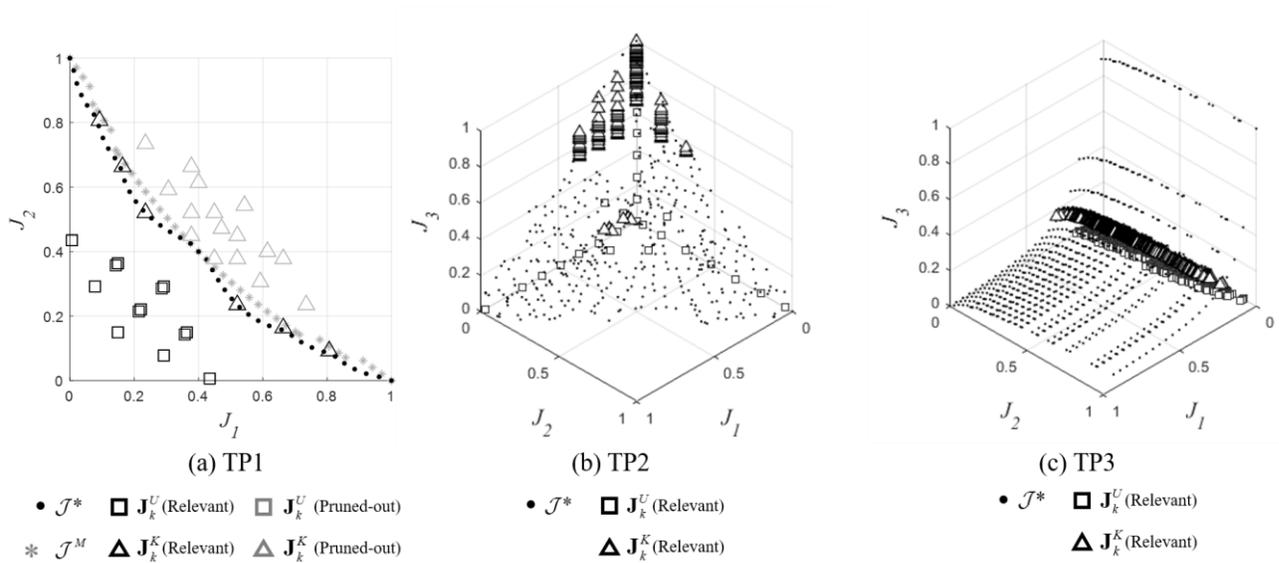

Figure 5: Pareto fronts of the test problems obtained by K-Pruning